\documentclass[11pt]{article}
\usepackage{amsmath,amssymb}

\def\eqn{Equation~}
\def\eqns{Equations~}

\begin{document}

%\markboth{William Dula and Ronald E. Mickens}{A note on the construction of a `valid'
%NSFD scheme for the Lotka-Volterra equations}

%%\articletype{NOTE}

\title{\bf A note on the construction of a `valid'
NSFD scheme for the Lotka-Volterra equations}

\author{William Dula$^{\rm a}$ and 
Ronald E. Mickens$^{\rm b}$\vspace{6pt}\thanks{Corresponding author.
Email: rohrs@math.gatech.edu}\\\vspace{6pt}
$^{\rm a}${Department of Mathematics, Morehouse College, 
Atlanta, GA 30314, USA};\\ $^{\rm b}${Department of Physics, Clark Atlanta
University, Atlanta, GA 30314, USA}}

\maketitle

\begin{abstract}
We demonstrate the construction of an explicit NSFD discretization for the
standard Lotka-Volterra equations.  In contrast to previous NSFD
schemes, our representation is dynamic consistent with respect to all the
essential properties of the differential equations and their solutions.
\bigskip


\end{abstract}

\section{Introduction}\label{sec1}
The main purpose of this note is to construct a ``valid" nonstandard
finite difference (NSFD) scheme for the Lotka-Volterra (L-V) equations 
\cite{6}
\begin{equation}\label{eq1}
\frac{dx}{dt}=ax-bxy,\quad \frac{dy}{dt}=-cy+dxy,
\end{equation}
where $(a,b,c,d)$ are non-negative parameters.  These equations provide
an elementary model for the interaction of a prey, $x$, with a
predator \cite{6,10}.

We use the word ``valid" to indicate that previous authors (see for example 
Roeger \cite{7,8}, who has worked extensively on NSFD schemes for L-V type
equations) have not used the full machinery of the NSFD methodology
to determine their particular discretizations of the L-V equations. 
A similar comment can be made for the recent work by Dang and Hoang
\cite{2}.  
However,
it must be emphasized that their obtained results are entirely correct;
they just did not make use of the maximum number  features based on
the NSFD scheme construction as formulated by Mickens \cite{3,4}.

It should be pointed out that an early publication to investigate the 
general Lotka-Volterra equations and their NSFD discretizations is AL-Kahby
et al.\ \cite{1}.  It is an outstanding work and covers the cases of
competitive, cooperative and predator-prey interactions.  In addition to
providing a broad range of possible NSFD schemes, it also gives rigorous
proofs of the dynamic consistency \cite{4} of the important features
of these schemes with respect to the differential equations. However, this
work was written prior to a full understanding as to how to 
construct the proper denominator functions (DF) for NSFD schemes
\cite{4}.  Consequently, while their DF's were correct, they were not
based on the (now known) NSFD methodology for constructing such
functions.  

In Section~\ref{sec2}, we give a brief overview of the important 
properties of the solutions to the L-V differential equations. 
Section~\ref{sec3} summarizes the difficulties which follow from using
a forward-Euler discretiztion for \eqns\eqref{eq1}.  Our main result is 
presented in Section~\ref{sec4}; here we construct a full, valid
NSFD scheme for the L-V equations.  Finally, in Section~\ref{sec5}, we 
summarize our main results and give several possible extensions of this
work.

\section{Properties of Lotka-Volterra equations}\label{sec2}
The following are several well known, basic features of the solutions
to the L-V differential equations \cite{6,10}:

(i) Positivity,
\begin{equation}\label{eq2}
\left.\begin{matrix}
x(0)=x_0 >0\\
y(0)=y_0 >0\end{matrix}\right\}\Longrightarrow
\begin{pmatrix}
x(t)>0\\ y(t)>0\end{pmatrix},\quad t>0.
\end{equation}

(ii) Two fixed-point (constant solutions) exist; they are located in the
$x-y$ phase-plane at
\begin{equation}\label{eq3}
FP_1=(0,0),\quad FP_2=\left(\frac cd,\frac ab\right).
\end{equation}

(iii) The first fixed-point is unstable and hyperbolic \cite{6,9},
while the second is a neutrally stable, center \cite{6,9}.

(iv) For $y(t)=0$, i.e., no predator is present, then 
\begin{equation}\label{eq4}
\frac{dx(t)}{dt} = ax(t), \quad x(0)=x_0>0,
\end{equation}
and $x(t)$ increases exponentially.  For $x(t)=0$, i.e., no prey is
available
\begin{equation}\label{eq5}
\frac{dy(t)}{dt} = -cy(t),\quad y(0)=y_0>0,
\end{equation}
then $y(t)$ decreases to zero, again exponentially.

(v) For small perturbations about the second fixed-point, i.e., 
\begin{subequations}\label{eq6}
\begin{align}
x(t) &= \frac cd + \alpha (t), \quad |\alpha (0)| \ll \frac cd,\label{eq6a}\\
y(t) &= \frac ab + \beta (t), \quad |\beta (0)| \ll \frac ab,\label{eq6b}
\end{align}
\end{subequations}
where $\alpha (t)$ and $\beta (t)$ satisfy (in the linear approximation)
the equation
\begin{equation}\label{eq7}
\frac{d^{\,2} z(t)}{dt^2}+(ac)z(t)=0,\quad
z(t)=\alpha (t)\mbox{ or } \beta (t),
\end{equation}
consequently, the motion is oscillatory.  In fact this result holds
for arbitrary initial conditions, $x_0>0$ and $y_0>0$.

\section{Forward-Euler scheme}\label{sec3}
A standard forward-Euler scheme for \eqns\eqref{eq1} is \cite{3,7}
\begin{equation}\label{eq8}
\frac{x_{k+1}-x_k}h =ax_k-bx_ky_k,\quad 
\frac{y_{k+1}-y_k}h =-cy_k+dx_ky_k,
\end{equation}
where $h=\Delta t$, $t\to hk=t_k$, $x(t)\to x_k$ and $y(t)\to y_k$, and
$k=(0,1,2,\dots)$.

Inspection of \eqn\eqref{eq8} shows that there exist initial conditions,
$x(0)=x_0>0$ and $y(0)=y_0>0$, and step-sizes, $h=\Delta t>0$, such that
the solutions, $(x_k,y_k)$, become negative.  Consequently, some of the five
conditions, given in Section~\ref{sec2} for the L-V differential equations,
may not be satisfied by the solutions of this forward-Euler 
discretization, and it must be concluded that this scheme is, in general,
not suitable or dynamic consistent \cite{5} with the original L-V
equations, and should not be used for their numerical integration.
It should also be noted that the previous work by Roeger \cite{7,8} had as
its focus the goal of eliminating these difficulties for various classes
of L-V differential equations.

Finally, an elementary calculation shows that for the scheme given by
\eqn\eqref{eq8}, the fixed-point  at $(c/d,a/b)$ is an unstable spiral
\cite{5,9}, i.e., near this fixed-point, solutions oscillate with
increasing values for the amplitude.  This result provides an additional
reason for the unsuitability of the forward-Euler discretization
of the L-V equations.

\section{Valid NSFD scheme}\label{sec4}
We now apply the full NSFD methodology \cite{3,4} to construct a 
discretization of \eqns\eqref{eq1}.

To calculate the two denominator functions \cite{3,4}, we need to
consider the respective linear terms in \eqns\eqref{eq1}, i.e.,
\begin{equation}\label{eq9}
\frac{dx}{dt} =ax,\quad \frac{dy}{dt} = -cy.
\end{equation}
The corresponding exact finite difference schemes are \cite{3}
\begin{equation}\label{eq10}
\frac{x_{k+1}-x_k}{\phi_1(a,h)} =ax_k,\quad
\frac{y_{k+1}-y_k}{\phi_2(c,h)} =-cy_k,
\end{equation}
where
\begin{equation}\label{eq11}
\phi_1(a,h) = \frac{e^{ah}-1}a,\quad 
\phi_2(c,h) = \frac{1-e^{-ch}}c,
\end{equation}
and $\phi_1(a,h)$ and $\phi_2(c,h)$ are the required denominator functions
needed for the discretizations of \eqns\eqref{eq1}.

The nonlinear term, $xy$, appears in both differential equations and must
be modeled nonlocally in the discrete representation.  To have positivity
for the $x$-variable, the term $xy$ takes the form
\begin{equation}\label{eq12}
xy\to x_{k+1}y_k.
\end{equation}
Since $xy$ must have the same structure where ever it appears, then it must
also take the form given by \eqn\eqref{eq12} in the $y$-variable equation.
Combining the results of \eqns\eqref{eq10} and \eqref{eq12}, it follows
that the NSFD schemes for the L-V equations are
\begin{subequations}\label{eq13}
\begin{align}
\frac{x_{k+1}-x_k}{\phi_1(a,h)} &= ax_k-bx_{k+1}y_k,\label{eq13a}\\
\frac{y_{k+1}-y_k}{\phi_2(c,h)} &= -cy_k+dx_{k+1}y_k,\label{eq13b}
\end{align}
\end{subequations}
where $\phi_1$ and $\phi_2$ are given in \eqn\eqref{eq11}.

Examination of \eqns\eqref{eq13} allows the following conclusions to be
reached:

(a) \eqns\eqref{eq13} have exactly the same two fixed-points as the
L-V differential equations, namely, 
\begin{equation}\label{eq14}
(\bar x, \bar y): (0,0)\mbox{ and } \left(\frac cd,\frac ab\right).
\end{equation}

(b) For $y_k\equiv 0$, then
\begin{equation}\label{eq15}
\frac{x_{k+1}-x_k}{\phi_1(a,h)} =ax_k\Longrightarrow x_k=
x_0 e^{ahk},
\end{equation}
and $x_k$  exponentially increases.

(c) For $x_k\equiv 0$, then
\begin{equation}\label{eq16}
\frac{y_{k+1}-y_k}{\phi_2(c,h)} =-cy_k\Longrightarrow y_k=
y_0 e^{-chk},
\end{equation}
and $y_k$  exponentially decreases to zero.

(d) Solving \eqns\eqref{eq13}, respectively, for $x_{k+1}$ and
$y_{k+1}$, gives
\begin{equation}\label{eq17}
x_{k+1} = \frac{e^{ah}x_k}{1+(b\phi_1)y_k},\quad
y_{k+1} = \left[ e^{-ch} + d\phi_2x_{k+1}\right]y_k.
\end{equation}
From the first relation, we note that if $x_k>0$ and $y_k>0$, then
$x_{k+1}>0$. Using this result in the second relation gives
$y_{k+1}>0$. Hence, the positivity condition holds.

(e) Now consider small perturbations about the fixed-point
$(\bar x,\bar y)=(c/d, a/b)$, i.e.,
\begin{equation}\label{eq18}
x_k=\bar x+\alpha_k,\quad y_k=\bar y+\beta_k,
\end{equation}
where 
\begin{equation}\label{eq19}
|\alpha_0|\ll \frac cd,\quad |\beta_0|\ll \frac ab.
\end{equation}
Substitution of \eqns\eqref{eq18} into \eqns\eqref{eq13} and retaining
only the linear terms, shows that $\alpha_k$ and $\beta_k$ both 
satisfy the second-order, linear, constant coefficient, difference
equation \cite{5}
\begin{equation}\label{eq20}
\frac{z_{k+1}-2z_k+z_{k-1}}{\phi_1(a,h)\phi_2(c,h)} + (ac)
z_k=0,
\end{equation}
where $z_k:\alpha_k$ or $\beta_k$.  Observe that this is a discrete
approximation for \eqn\eqref{eq7}.  Hence, we conclude that for small
perturbations about the nontrivial fixed-point, the trajectory
in the $(x_k,y_k)$ phase-space is periodic.  (To obtain this result in detail,
see the technique presented in Mickens \cite{5}.)

Our conclusion is that the NSFD discretization of the L-V \eqns\eqref{eq1}
is fully dynamic consistent with respect to all five conditions specified
in Section~\ref{sec2}.

One of the major differences between the work of Roeger \cite{7,8}
and the current results is that Roeger uses the simple
expression, $h$, for the denominator functions; we have explicitly 
calculated the required functions.

\section{Extensions}\label{sec5}
We have in this note constructed an improved NSFD scheme for the L-V
predator-prey differential equations.  Further, we have shown it to be
dynamic consistent with all of the important features of the solutions 
to the standard L-V equations.  
While we do not give the calculations here, it is methodology easy to
extend the results
to the generalized L-V equations investigated by Roeger \cite{7,8}
and AL-Kahby et al.\ \cite{1}.  Finally, another important case to
investigate is where fractional-order derivatives appear \cite{1}.

\end{document}